\documentclass[11pt]{amsart}
\usepackage{amsmath,amssymb,amsfonts,url}

\numberwithin{equation}{section}

\newtheorem{thm}{Theorem}[section]
\newtheorem{lem}{Lemma}[section]
\newtheorem{rem}{Remark}[section]

\begin{document}
\title[Gradient estimate]{Local Gradient Estimate for $p$-harmonic functions
on Riemannian Manifolds}
\author{Xiaodong Wang}
\address{Department of Mathematics\\
Michigan State University\\
East Lansing, MI 48824}
\email{xwang@math.msu.edu}
\author{Lei Zhang}
\address{Department of Mathematics\\
University of Florida\\
358 Little Hall, P. O. Box 118105\\
Gainesville, FL 32611-8105}
\email{leizhang@ufl.edu}
\thanks{Wang acknowledges support from NSF grant DMS-0905904. Zhang is
supported in part by NSF Grant 0900864(1027628)}
\date{\today }

\begin{abstract}
For positive $p$-harmonic functions on Riemannian manifolds, we derive a
gradient estimate and Harnack inequality with constants depending only on
the lower bound of the Ricci curvature, the dimension $n$, $p$ and the
radius of the ball on which the function is defined. Our approach is based
on a careful application of the Moser iteration technique and is different
from Cheng-Yau's method \cite{CY} employed by Kostchwar and Ni \cite{kn}, in
which a gradient estimate for positive $p$-harmonic functions is derived
under the assumption that the sectional curvature is bounded from below.
\end{abstract}

\maketitle

\section{Introduction}

\bigskip The study of harmonic functions on Riemannian manifolds has been
one of the central subjects in geometric analysis. In their classical work
Cheng-Yau \cite{SY} derived the following gradient estimate for positive
harmonic functions on Riemannian manifolds:

\emph{Theorem A (Cheng-Yau) Let }$M$ \emph{be an }$n$\emph{-dimensional
complete Riemannian manifold with }$Ric\geq -\left( n-1\right) \kappa $\emph{%
, where }$\kappa \geq 0$ \emph{is a constant. Suppose that }$u$ \emph{is a
positive harmonic function on a geodesic ball }$B\left( o,R\right) $\emph{.
Then}%
\begin{equation}
\sup_{B\left( o,R/2\right) }\frac{\left\vert \nabla u\right\vert }{u}\leq
C_{n}\frac{1+R\sqrt{\kappa }}{R},  \label{yau2}
\end{equation}%
\ \emph{where }$C_{n}$\emph{\ is a constant depending only on }$n$\emph{.}

\bigskip

An important feature of Cheng-Yau's estimate is that the RHS (which stands
for the right hand side) of (\ref{yau2}) depends only on $n$, $k$ and $R$,
it does not depend on the lower bound of the injectivity radius or a global
coordinate system. From PDE viewpoints, deriving a Harnack inequality
requires some bounds on the coefficients in some fixed coordinate system,
thus not suitable for many problems defined on manifolds. We also observe
that the RHS of (\ref{yau2}) is optimal in the sense that even for $\kappa>0$%
, the bound stays bounded when $R\rightarrow \infty$.

There are two major ingredients in the proof of Theorem A. First a Bochner
formula is used to derive a lower bound of the Laplacian of $|\nabla u|^{2}$
for a harmonic function $u$ in terms of the lower bound of the Ricci tensor.
The second major ingredient is a clever application of the maximum
principle. The trick is to multiply $|\nabla u|^{2}$ by a cut-off function,
derive a new differential inequality for the product and then apply the
maximum principle. The cut-off function is constructed using the distance
function. As a result, the new differential inequality involves the
Laplacian of the distance function. As is well known, the Riemannian
distance function is uniformly Lipschitz and its Laplacian has an upper
bound depending on the lower bound of the Ricci tensor.

Cheng-Yau's approach turned out to be very useful and some important results
for other problems are deeply influenced by Theorem A. For example, P. Li 
\cite{pli} obtained the sharp lower bound for the first eigenvalue of a
manifold, which was later generalized by Li-Yau \cite{liyau1}. Similar
results were also obtained by Li-Yau \cite{liyau2} for heat equations. S. Y.
Cheng \cite{cheng} and H. I. Choi \cite{choi} obtained gradient estimates
for harmonic mappings, etc. We refer to \cite{SY, schoen} and the more
recent survey \cite{Li} for an overview of the subject.

$p$-harmonic functions are natural extensions of harmonic functions from a
variational point of view. It has been extensively studied because of its
various interesting features and applications. Compared with the theory for
harmonic functions the study of $p$-harmonic functions is generally harder
because the equation, even though elliptic, is degenerate and the regularity
results are far weaker (see, for example \cite{tolk}). Recently, there has
been renewed interest in $p$-harmonic functions. In particular R. Moser \cite%
{M} established a nice connection between $p$-harmonic functions and the
inverse mean curvature flow. In a recent paper \cite{kn} Kotschwar and Ni
derived, among other things, a local gradient estimate for $p$-harmonic
functions under the assumption that the sectional curvature is bounded from
below. It is remarkable that the constant in their estimate does not blowup
when $p\rightarrow 1$, which leads to interesting results on the inverse
mean curvature flow problems. Their proof follows the same strategy
introduced by Cheng-Yau \cite{CY} for harmoinc functions (i.e. $p=2$).
However for general $p$-harmonic functions, the computation involves the
full Hessian of the distance function when the cut-off function is
introduced. As a result, a lower bound on the sectional curvature has to be
assumed in \cite{kn}.

Kotschwar and Ni speculated that their estimate may hold if only a lower
bound on the Ricci tensor is assumed. The main result of this paper is to
establish the following theorem:

\begin{thm}
\label{p-harm} Let $\left( M^{n},g\right) $ be a complete Riemannian
manifold with $\mathrm{Ric}\geq -\left( n-1\right) \kappa $.\bigskip\ Assume
that $v$ is a positive $p$-harmonic function on the ball $B(o,R)\subset M$.
Then there exists a constant $C_{p,n}$ such that%
\begin{equation*}
\frac{|\nabla v|}{v}\leq C_{p,n}\left( 1+\sqrt{\kappa }R\right) /R\quad %
\mbox{on}\quad B(o,R/2).
\end{equation*}
\end{thm}

The proof of Theorem \ref{p-harm} will be presented in section two. As far
as the second major ingredient of Cheng-Yau's proof is concerned, our
approach follows a different strategy by carefully using the Moser iteration
technique. This approach only involves differentiating the distance function
once and hence bypasses the difficulty of handling the full Hessian of the
distance function. In the special case $p=2$, when $p$-harmonic function are
just harmonic functions, Theorem \ref{p-harm} is exactly Cheng-Yau's
theorem. An immediate consequence of Theorem \ref{p-harm} is the following
Harnack inequality.

\begin{thm}
\bigskip Let $\left( M^{n},g\right) $ be a complete Riemannian manifold with 
$\mathrm{Ric}\geq -\left( n-1\right) \kappa $.\bigskip\ Assume that $v$ is a
positive $p$-harmonic function on the ball $B(o,R)\subset M$. Then there
exists a constant $C_{p,n}$ such that for any $x,y\in B(o,R/2)$,%
\begin{equation*}
v\left( x\right) /v\left( y\right) \leq e^{C_{p,n}\left( 1+\sqrt{\kappa }%
R\right) }.
\end{equation*}
\end{thm}

\bigskip

It follows that if $\mathrm{Ric}\geq 0$, then we have a uniform constant $%
c_{p,n}$ (independent of $R$) s.t. 
\begin{equation}
\sup_{B(o,R/2)}v\leq c_{p,n}\inf_{B(o,R/2)}v.  \label{harn}
\end{equation}%
This was already proved by Rigoli, Salvatori, and Vignati \cite{RSV}. In
fact, they proved the stronger result that (\ref{harn}) holds provided that
the volume is doubling and a weak Poincare inequality holds. See also \cite%
{Hol}.

Another standard application of Theorem \ref{p-harm} is the following
Liouville theorem, which was also deduced from the Harnack inequality in 
\cite{RSV}.

\emph{\ Let }$u$\emph{\ be a }$p$\emph{-harmonic function bounded from above
or below on a complete Riemannian manifold with non-negative Ricci tensor,
then }$u$\emph{\ is constant. } \medskip

Finally we point out that our constant $C_{p,n}$ in Theorem \ref{p-harm}
becomes unbounded as $p\rightarrow 1$, while in Kotschwar-Ni's result, all
the constants stay bounded when $p\rightarrow 1$. We do not know if the
method can be tweaked to remove this defect.

\section{The gradient estimate for the $p$-harmonic functions}

$p$-harmonic functions arise naturally as critical points of the $%
L^{p}\,\,(p>1)$ norm of the gradient. Let $\left( M^{n},g\right) $ be a
complete Riemannian manifold and $\Omega \subset M$ an open set. A function $%
v\in W_{loc}^{1,p}\left( \Omega \right) $ is $p$-harmonic if%
\begin{equation*}
\mathrm{div}\left( \left\vert \nabla v\right\vert ^{p-2}\nabla v\right) =0
\end{equation*}%
in the weak sense, i.e.%
\begin{equation*}
\int_{U}\left\vert \nabla v\right\vert ^{p-2}\left\langle \nabla v,\nabla
\xi \right\rangle =0
\end{equation*}%
for all $\xi \in W_{0}^{1,p}\left( \Omega \right) $. By \cite{tolk} for
example, $v$ must be $C^{1,\alpha }$. Moreover $v\in W_{loc}^{2,2}$ if $%
p\geq 2$; $v\in W_{loc}^{2,p}$ if $1<p<2$. Away from $\left\{ \nabla
v=0\right\} $, $v$ is in fact smooth.

Suppose that $v$ is positive. Set $u=-\left( p-1\right) \log v$. Then $u$
satisfies%
\begin{equation}
\mathrm{div}\left( \left\vert \nabla u\right\vert ^{p-2}\nabla u\right)
=\left\vert \nabla u\right\vert ^{p}.  \label{equ}
\end{equation}%
Let $f=\left\vert \nabla u\right\vert ^{2}$. \ We define%
\begin{equation*}
\mathcal{L}\left( \psi \right) =\mathrm{div}\left( f^{p/2-1}A\left( \nabla
\psi \right) \right) -pf^{p/2-1}\left\langle \nabla u,\nabla \psi
\right\rangle ,
\end{equation*}%
where 
\begin{equation*}
A=id+\left( p-2\right) \frac{\nabla u\otimes \nabla u}{\left\vert \nabla
u\right\vert ^{2}}.
\end{equation*}

We need the following lemma from \cite{kn} and the proof is by direct
calculation.

\begin{lem}
\label{nilem} 
\begin{equation*}
\mathcal{L}\left( f\right) =2f^{p/2-1}\left( \left\vert D^{2}u\right\vert
^{2}+\mathrm{Ric}\left( \nabla u,\nabla u\right) \right) +\left( \frac {p}{2}%
-1\right) \left\vert \nabla f\right\vert ^{2}f^{p/2-2}.
\end{equation*}
\end{lem}

\begin{rem}
Lemma \ref{nilem} holds point-wisely in $\{x:\,\,f(x)>0\}$. From the
gradient estimate of \cite{tolk} we know that $f=|\nabla u|^{2}\in C^{\alpha
}$ for some $\alpha >0$ and $f\in W_{loc}^{1,\beta }$ for some $\beta >1$.
\end{rem}

We choose a local orthonormal frame $\left\{ e_{i}\right\} $ with $%
e_{1}=\nabla u/\left\vert \nabla u\right\vert $. Then (\ref{equ}) takes the
following form%
\begin{equation*}
\left( p-1\right) u_{11}+\sum_{i=2}^{n}u_{ii}=f.
\end{equation*}%
Therefore%
\begin{align*}
\left\vert D^{2}u\right\vert ^{2}& \geq
u_{11}^{2}+2\sum_{i=2}^{n}u_{1i}^{2}+\sum_{i=2}^{n}u_{ii}^{2} \\
& \geq u_{11}^{2}+2\sum_{i=2}^{n}u_{1i}^{2}+\frac{1}{n-1}\left(
\sum_{i=2}^{n}u_{ii}\right) ^{2} \\
& =u_{11}^{2}+2\sum_{i=2}^{n}u_{1i}^{2}+\frac{1}{n-1}\left( f-\left(
p-1\right) u_{11}\right) ^{2} \\
& =\frac{1}{n-1}f^{2}-\frac{2\left( p-1\right) }{n-1}fu_{11}+\left( 1+\frac{%
\left( p-1\right) ^{2}}{n-1}\right) u_{11}^{2}+2\sum_{i=2}^{n}u_{1i}^{2} \\
& \geq \frac{1}{n-1}f^{2}-\frac{2\left( p-1\right) }{n-1}fu_{11}+a_{0}%
\sum_{i=1}^{n}u_{1i}^{2},
\end{align*}%
where $a_{0}=1+\min \left( \frac{\left( p-1\right) ^{2}}{n-1},1\right) >1$.
Using the identities 
\begin{align*}
2fu_{11}& =\left\langle \nabla u,\nabla f\right\rangle , \\
\sum_{i=1}^{n}u_{1i}^{2}& =\frac{1}{4}\frac{\left\vert \nabla f\right\vert
^{2}}{f},
\end{align*}%
we end up with 
\begin{equation*}
\left\vert D^{2}u\right\vert ^{2}\geq \frac{1}{n-1}f^{2}-\frac{\left(
p-1\right) }{n-1}\left\langle \nabla u,\nabla f\right\rangle +\frac{a_{0}}{4}%
\frac{\left\vert \nabla f\right\vert ^{2}}{f}.
\end{equation*}%
Assume that $\mathrm{Ric}\geq -\left( n-1\right) \kappa $. Therefore 
\begin{align}
\mathcal{L}\left( f\right) & \geq -2\left( n-1\right) \kappa f^{p/2}+\left( 
\frac{p+a_{0}}{2}-1\right) \left\vert \nabla f\right\vert ^{2}f^{p/2-2}
\label{930e1} \\
& +\frac{2}{n-1}f^{p/2+1}-\frac{2\left( p-1\right) }{n-1}f^{p/2-1}\left%
\langle \nabla u,\nabla f\right\rangle  \notag \\
& \geq -2\left( n-1\right) kf^{p/2}+\frac{2}{n-1}f^{p/2+1}-\frac{2\left(
p-1\right) }{n-1}f^{p/2-1}\left\langle \nabla u,\nabla f\right\rangle  \notag
\end{align}

Equation (\ref{930e1}) holds wherever $f$ is strictly positive. Let $%
K=\{x\in \Omega :f(x)=0\}$. Then for any nonnegative function $\psi $ with
compact support in $\Omega \setminus K$, we have 
\begin{eqnarray}
&&\int_{\Omega }\left\langle f^{p/2-1}\nabla f+(p-2)f^{p/2-2}\left\langle
\nabla u,\nabla f\right\rangle \nabla u,\nabla \psi \right\rangle
\label{lfweak} \\
&&+p\int_{\Omega }f^{p/2-1}\left\langle \nabla u,\nabla f\right\rangle \psi +%
\frac{2}{n-1}\int_{\Omega }f^{p/2+1}\psi  \notag \\
&\leq &2(n-1)k\int_{\Omega }f^{p/2}\psi +\frac{2(p-1)}{n-1}\int_{\Omega
}f^{p/2-1}\left\langle \nabla u,\nabla f\right\rangle \psi .  \notag
\end{eqnarray}%
In particular, let $\epsilon >0$ and $\psi =f_{\epsilon }^{b}\eta ^{2}$
where $f_{\epsilon }=(f-\epsilon )^{+}$, $\eta \in C_{0}^{\infty }(B_{R})$
is non-negative, $b>1$ is to be determined later. Then direct computation
yields 
\begin{equation*}
\nabla \psi =bf_{\epsilon }^{b-1}\nabla f\eta ^{2}+2f_{\epsilon }^{b}\eta
\nabla \eta .
\end{equation*}%
Using the above in (\ref{lfweak}) we obtain 
\begin{eqnarray}
&&b\int_{B_{R}}\bigg (f^{p/2-1}f_{\epsilon }^{b-1}|\nabla
f|^{2}+(p-2)f^{p/2-2}\left\langle \nabla u,\nabla f\right\rangle
^{2}f_{\epsilon }^{b-1}\bigg )\eta ^{2}  \label{930e3} \\
&+&2(p-2)\int_{B_{R}}f^{p/2-2}\left\langle \nabla u,\nabla f\right\rangle
f_{\epsilon }^{b}\eta (\nabla u,\nabla \eta )+2\int_{\Omega
}f^{p/2-1}f_{\epsilon }^{b}\eta \left\langle \nabla f,\nabla \eta
\right\rangle  \notag \\
&+&p\int_{\Omega }f^{p/2-1}\left\langle \nabla u,\nabla f\right\rangle
f_{\epsilon }^{b}\eta ^{2}+\frac{2}{n-1}\int_{\Omega }f^{p/2+1}f_{\epsilon
}^{b}\eta ^{2}  \notag \\
&\leq &2(n-1)\kappa \int_{\Omega }f^{p/2-1}f_{\epsilon }^{b}\eta ^{2}+\frac{%
2(p-1)}{n-1}\int_{\Omega }f^{p/2-1}\left\langle \nabla u,\nabla
f\right\rangle f_{\epsilon }^{b}\eta ^{2}.  \notag
\end{eqnarray}%
Since $u\in C^{1,\alpha }$, $f\in C^{\alpha }$ and $\nabla f\in L^{\beta
}(\Omega )$ for some $\alpha >0$ and $\beta >1$, we see that except for the
first term, all the other terms converge to the corresponding form without $%
\epsilon $. For the first term, observe that 
\begin{equation*}
f^{p/2-1}f_{\epsilon }^{b-1}|\nabla f|^{2}+(p-2)f^{p/2-2}\left\langle \nabla
u,\nabla f\right\rangle ^{2}f_{\epsilon }^{b-1}\geq
a_{1}f^{p/2-1}f_{\epsilon }^{b-1}|\nabla f|^{2}
\end{equation*}%
where $a_{1}=1$ if $p\geq 2$ and $a_{1}=(p-1)$ if $p\in (1,2)$. Thus by
passing $\epsilon $ to $0$ we have 
\begin{eqnarray}
&&ba_{1}\int_{\Omega }f^{p/2+b-2}|\nabla f|^{2}\eta ^{2}  \label{930e2} \\
&+&2(p-2)\int_{\Omega }f^{\frac{p-4}{2}+b}\left\langle \nabla u,\nabla
f\right\rangle \eta \left\langle \nabla u,\nabla \eta \right\rangle
+2\int_{\Omega }f^{\frac{p-2}{2}+b}\eta \left\langle \nabla f,\nabla \eta
\right\rangle  \notag \\
&+&p\int_{\Omega }f^{\frac{p-2}{2}+b}\left\langle \nabla u,\nabla
f\right\rangle \eta ^{2}+\frac{2}{n-1}\int_{\Omega }f^{\frac{p+2}{2}+b}\eta
^{2}  \notag \\
&\leq &2(n-1)k\int_{\Omega }f^{\frac{p-2}{2}+b}\eta ^{2}+\frac{2(p-1)}{n-1}%
\int_{\Omega }f^{\frac{p-2}{2}+b}\left\langle \nabla u,\nabla f\right\rangle
\eta ^{2}.  \notag
\end{eqnarray}

From now on we use $a_{1},a_{2},\cdots $ etc. to denote constants depending
only on $p$ and $n$. Combining terms in (\ref{930e2}) using the definition
of $f$ we have 
\begin{eqnarray}
&&a_{1}b\int_{\Omega }f^{p/2+b-2}|\nabla f|^{2}\eta ^{2}+\frac{2}{n-1}%
\int_{\Omega }f^{p/2+1+b}\eta ^{2}  \label{921e1} \\
&\leq &2(n-1)\kappa \int_{\Omega }f^{p/2+b}\eta ^{2}+a_{2}\int_{\Omega }f^{%
\frac{p-1}{2}+b}|\nabla f|\eta ^{2}  \notag \\
&&+a_{3}\int_{\Omega }f^{p/2+b-1}|\nabla f||\nabla \eta |\eta .  \notag
\end{eqnarray}%
For $R_{3}$ (the third term on the RHS, $L_{1},L_{2},R_{1}$ etc. are
understood similarly) in (\ref{921e1}) we have 
\begin{equation*}
|R_{3}|\leq \frac{a_{1}b}{4}\int_{\Omega }f^{b+p/2-2}|\nabla f|^{2}\eta ^{2}+%
\frac{a_{4}}{b}\int_{\Omega }|\nabla \eta |^{2}f^{b+p/2}.
\end{equation*}%
Also by Cauchy's inequality $R_{2}$ of (\ref{921e1}) can be estimated as 
\begin{equation*}
|R_{2}|\leq \frac{a_{1}b}{4}\int_{\Omega }f^{p/2+b-2}|\nabla f|^{2}\eta ^{2}+%
\frac{a_{5}}{b}\int_{\Omega }f^{b+p/2+1}\eta ^{2}.
\end{equation*}%
With these two inequalities we obtain 
\begin{eqnarray}
&&\frac{a_{1}b}{2}\int_{\Omega }f^{p/2+b-2}|\nabla f|^{2}\eta ^{2}+\frac{2}{%
n-1}\int_{\Omega }f^{p/2+1+b}\eta ^{2}  \label{923e1} \\
&\leq &\int_{\Omega }\left( 2(n-1)\kappa \eta ^{2}+\frac{a_{4}}{b}|\nabla
\eta |^{2}\right) f^{p/2+b}+\frac{a_{5}}{b}\int_{\Omega }f^{b+p/2+1}\eta
^{2}.  \notag
\end{eqnarray}%
By requiring 
\begin{equation}
\frac{a_{5}}{b}<\frac{1}{n-1}  \label{assumb}
\end{equation}%
we see that the last term in on the RHS of (\ref{923e1}) is majorized by the
last term on the LHS. Therefore we have 
\begin{eqnarray*}
&&\frac{a_{1}b}{2}\int_{\Omega }f^{p/2+b-2}|\nabla f|^{2}\eta ^{2}+\frac{1}{%
n-1}\int_{\Omega }f^{p/2+1+b}\eta ^{2} \\
&\leq &2(n-1)\kappa \int_{\Omega }f^{p/2+b}\eta ^{2}+\frac{a_{4}}{b}%
\int_{\Omega }|\nabla \eta |^{2}f^{p/2+b}.
\end{eqnarray*}%
For the first term on the LHS we use 
\begin{equation*}
|\nabla (f^{p/4+b/2}\eta )|^{2}\leq \frac{1}{2}(\frac{p}{2}%
+b)^{2}f^{p/2+b-2}|\nabla f|^{2}\eta ^{2}+2f^{p/2+b}|\nabla \eta |^{2}.
\end{equation*}%
From the above we obtain 
\begin{equation}
\int_{\Omega }|\nabla (f^{p/4+b/2}\eta )|^{2}+d_{1}\int_{\Omega
}f^{p/2+1+b}\eta ^{2}\leq \kappa d_{2}\int_{\Omega }f^{p/2+b}\eta
^{2}+a_{7}\int_{\Omega }|\nabla \eta |^{2}f^{b+p/2}.  \label{929e1}
\end{equation}%
where $d_{1}\sim b,d_{2}\sim b$ (recall $b>1$, $d_{1}\sim b$ means $d_{1}$
is comparable to $b$, $d_{2}\sim b$ is understood the same way).

The following Sobolev embedding theorem of Saloff-Coste plays an important
role in our approach:

\medskip \emph{Theorem B (Theorem 3.1 of \cite{saloff}) Let }$\left(
M^{n},g\right) $\emph{\ be a complete Riemannian manifold with }$Ric\geq
-\left( n-1\right) \kappa $\emph{. For }$n>2$\emph{, there exists }$C$\emph{%
, depending only on }$n$\emph{, such that for all }$B\subset M$\emph{\ of
radius }$R$\emph{\ and volume }$V$\emph{\ we have for }$f\in C_{0}^{\infty
}\left( B\right) $%
\begin{equation*}
\left( \int \left\vert f\right\vert ^{2q}\right) ^{1/q}\leq e^{C\left( 1+%
\sqrt{\kappa }R\right) }V^{-2/n}R^{2}\left( \int \left\vert \nabla
f\right\vert ^{2}+R^{-2}f^{2}\right) ,
\end{equation*}%
\emph{\ where }$q=n/(n-2)$\emph{. For }$n=2$\emph{, the above inequality
holds with }$n$\emph{\ replaced by any fixed }$n^{\prime }>2$\emph{. }

\medskip

From now on, we assume $\Omega =B\left( o,R\right) $. Theorem B gives 
\begin{eqnarray}
&&\left( \int_{\Omega }f^{\frac{n(p/2+b)}{n-2}}\eta ^{\frac{2n}{n-2}}\right)
^{\left( n-2\right) /n}  \label{929e2} \\
&\leq &e^{c_{0}(1+\sqrt{\kappa }R)}V^{-\frac{2}{n}}\bigg (R^{2}\int_{\Omega
}|\nabla (f^{p/4+b/2}\eta )|^{2}+\int_{\Omega }f^{p/2+b}\eta ^{2}\bigg ). 
\notag
\end{eqnarray}%
where $c_{0}(n,p)>0$ depends only on $n,p$. Let $b_{0}=c_{1}(n,p)(1+\sqrt{%
\kappa }R)$ with $c_{1}(n,p)$ large enough to make $b_{0}$ satisfy (\ref%
{assumb}), then (\ref{929e1}) and (\ref{929e2}) combined gives 
\begin{eqnarray}
&&\left( \int_{\Omega }f^{\frac{n(p/2+b)}{n-2}}\eta ^{\frac{2n}{n-2}}\right)
^{\left( n-2\right) /n}+a_{8}bR^{2}e^{c_{2}b_{0}}V^{-2/n}\int_{\Omega
}f^{p/2+1+b}\eta ^{2}  \label{921moser} \\
&\leq &a_{9}b_{0}^{2}be^{c_{2}b_{0}}V^{-2/n}\int_{\Omega }f^{p/2+b}\eta
^{2}+a_{10}e^{c_{2}b_{0}}V^{-2/n}R^{2}\int_{\Omega }|\nabla \eta
|^{2}f^{p/2+b}  \notag
\end{eqnarray}

\begin{lem}
\label{923lem1} Let $b_{1}=(b_{0}+\frac{p}{2})\frac{n}{n-2}$. Then there
exists $c_{3}(n,p)>0$ such that 
\begin{equation}
\Vert f\Vert _{L^{b_{1}}(B_{3R/4})}\leq c_{3}\frac{b_{0}^{2}}{R^{2}}%
V^{1/b_{1}}.  \label{923e2}
\end{equation}
\end{lem}

\noindent \textbf{Proof of Lemma \ref{923lem1}:} Let $b=b_{0}$ in (\ref%
{921moser}), then by comparing $L_{2}$ and $R_{1}$ of (\ref{921moser}) we
observe that 
\begin{equation*}
a_{9}b_{0}^{3}f^{p/2+b_{0}}<\frac{1}{2}a_{8}b_{0}R^{2}f^{p/2+1+b_{0}}
\end{equation*}%
if $f>a_{11}b_{0}^{2}R^{-2}$. Thus in the evaluation of $R_{1}$ we decompose 
$\Omega $ into two subregions, one over the places where $f\leq
a_{11}b_{0}^{2}R^{-2}$ and the second region is the complement of the first
region. With this decomposition we have 
\begin{equation*}
R_{1}\leq a_{12}^{b_{0}}b_{0}^{3}\left( \frac{b_{0}}{R}\right)
^{p+2b_{0}}e^{c_{2}b_{0}}V^{1-2/n}+\frac{L_{2}}{2}.
\end{equation*}%
Now (\ref{921moser}) with $b=b_{0}$ can be written as 
\begin{eqnarray}
&&\left( \int_{\Omega }f^{\frac{n(p/2+b_{0})}{n-2}}\eta ^{\frac{2n}{n-2}%
}\right) ^{\left( n-2\right) /n}+\frac{a_{8}}{2}%
b_{0}R^{2}e^{c_{2}b_{0}}V^{-2/n}\int_{\Omega }f^{p/2+1+b_{0}}\eta ^{2}
\label{923e3} \\
&\leq &a_{12}^{b_{0}}b_{0}^{3}\left( \frac{b_{0}}{R}\right)
^{p+2b_{0}}e^{c_{2}b_{0}}V^{1-2/n}+a_{10}e^{c_{2}b_{0}}V^{-2/n}R^{2}\int_{%
\Omega }|\nabla \eta |^{2}f^{p/2+b_{0}}  \notag
\end{eqnarray}%
Now we choose $\eta $ to make $R_{2}$ in (\ref{923e3}) dominated by the LHS.
Let $\eta _{1}\in C_{0}^{\infty }(B_{R})$ satisfy 
\begin{equation*}
0\leq \eta _{1}\leq 1,\quad \eta _{1}\equiv 1\quad \mbox{in}\quad
B_{3R/4},\quad |\nabla \eta _{1}|\leq C(n)/R.
\end{equation*}%
Let $\eta =\eta _{1}^{m}$ where $m=b_{0}+\frac{p}{2}+1$. Direct computation
shows 
\begin{equation}
R^{2}|\nabla \eta |^{2}\leq a_{13}b_{0}^{2}\,\eta ^{\frac{2b_{0}+p}{%
b_{0}+p/2+1}}.  \label{923e4}
\end{equation}%
By (\ref{923e4}) and Young's inequality, the $R_{2}$ of (\ref{923e3}) can be
written as 
\begin{eqnarray*}
&&a_{10}\,R\,^{2}\int |\nabla \eta |^{2}f^{b_{0}+p/2} \\
&\leq &a_{14}\,b_{0}^{2}\,\int_{\Omega }f^{b_{0}+p/2}\eta ^{\frac{2b_{0}+p}{%
b_{0}+p/2+1}} \\
&\leq &a_{14}\,b_{0}^{2}\left( \int_{\Omega }f^{b_{0}+p/2+1}\eta ^{2}\right)
^{\frac{b_{0}+p/2}{b_{0}+p/2+1}}V^{\frac{1}{b_{0}+p/2+1}} \\
&\leq &\frac{a_{8}b_{0}}{2}\,R^{2}\,\int_{\Omega }f^{b_{0}+p/2+1}\eta
^{2}+a_{15}^{b_{0}+p/2}\frac{b_{0}^{b_{0}+p/2+2}}{R^{2b_{0}+p}}V.
\end{eqnarray*}

With the estimates on $R_{1},R_{2}$ we arrive at 
\begin{equation}
\bigg (\int_{B_{3R/4}}f^{(b_{0}+p/2)n/\left( n-2\right) }\bigg )^{\left(
n-2\right) /n}\leq a_{16}^{b_{0}}e^{c_{2}b_{0}}V^{1-2/n}b_{0}^{3}\left( 
\frac{b_{0}}{R}\right) ^{p+2b_{0}}.  \label{923e7}
\end{equation}

Recall $b_{1}=(b_{0}+\frac{p}{2})\frac{n}{n-2}$. Taking the $1/(b_{0}+p)$
root on both sides of (\ref{923e7}) we have 
\begin{equation*}
\Vert f\Vert _{L^{b_{1}}(B_{3R/4})}\leq a_{17}V^{\frac{1}{b_{1}}%
}b_{0}^{2}/R^{2}.
\end{equation*}%
Lemma \ref{923lem1} is established. $\Box $

\bigskip

Now we go back to (\ref{921moser}), by ignoring $L_{2}$ we have

\begin{eqnarray}
&&\left( \int_{\Omega }f^{\frac{n(p/2+b)}{n-2}}\eta ^{\frac{2n}{n-2}}\right)
^{\left( n-2\right) /n}  \label{923e9} \\
&\leq &a_{17}\frac{e^{c_{2}b_{0}}}{V^{2/n}}\int_{\Omega }\bigg (%
(b_{0}^{2}b\eta ^{2}+R^{2}|\nabla \eta |^{2}\bigg )f^{b+p/2}.  \notag
\end{eqnarray}%
To apply the Moser iteration we set 
\begin{equation*}
b_{l+1}=b_{l}\frac{n}{n-2},\quad \Omega _{l}=B(o,\frac{R}{2}+\frac{R}{4^{l}}%
),\quad l=1,2..
\end{equation*}%
and choose $\eta _{l}\in C_{0}^{\infty }\left( \Omega \right) $ s.t. 
\begin{equation*}
\eta _{l}\equiv 1\,\,\mbox{in}\,\,\Omega _{l+1},\quad \eta _{l}\equiv 0\,\,%
\mbox{in}\,\,\Omega \setminus \Omega _{l},\quad |\nabla \eta _{l}|\leq \frac{%
C4^{l}}{R},\quad 0\leq \eta _{l}\leq 1.
\end{equation*}%
Then in (\ref{923e9}), by letting $b+\frac{p}{2}=b_{l}$, $\eta =\eta _{l}$
we have 
\begin{equation*}
\left( \int_{\Omega _{l+1}}f^{b_{l+1}}\right) ^{\frac{1}{b_{l+1}}}\leq
\left( a_{17}\frac{e^{c_{2}b_{0}}}{V^{\frac{2}{n}}}\right) ^{\frac{1}{b_{l}}%
}\left( \int_{\Omega _{l}}\left( b_{0}^{2}b_{l}+R^{2}|\nabla \eta
_{l}|^{2}\right) f^{b_{l}}\right) ^{\frac{1}{b_{l}}}.
\end{equation*}%
By the estimate of $|\nabla \eta _{l}|$ 
\begin{equation}
\Vert f\Vert _{L^{b_{l+1}}(\Omega _{l+1})}\leq \left( a_{17}\frac{%
e^{c_{2}b_{0}}}{V^{2/n}}\right) ^{1/b_{l}}\left(
b_{0}^{2}b_{l}+16^{l}\right) ^{1/b_{l}}\Vert f\Vert _{L^{b_{l}}(\Omega
_{l})}.  \label{923e10}
\end{equation}%
Notice that $\sum_{l=1}^{\infty }\frac{1}{b_{l}}=\frac{n}{2b_{1}}$, then (%
\ref{923e10}) leads to 
\begin{eqnarray}
&&\Vert f\Vert _{L^{\infty }(B_{R/2})}  \label{923e11} \\
&\leq &\left( a_{18}\frac{e^{c_{2}b_{0}}}{V^{2/n}}\right)
^{\sum_{l=1}^{\infty }1/b_{l}}\prod\limits_{l=1}^{\infty }\left(
b_{0}^{3}\left( \frac{n}{n-2}\right) ^{l}+16^{l}\right) ^{1/b_{l}}\Vert
f\Vert _{L^{b_{1}}(B_{3R/4})}  \notag \\
&\leq &a_{19}\frac{e^{\frac{nc_{2}b_{0}}{2b_{1}}}}{V^{1/b_{1}}}b_{0}^{\frac{%
3n}{2b_{1}}}\Vert f\Vert _{L^{b_{1}}(B_{3R/4})}.  \notag
\end{eqnarray}%
Using Lemma \ref{923lem1} in (\ref{923e11}) we obtain 
\begin{equation}
\Vert f\Vert _{L^{\infty }(B_{R/2})}\leq a_{20}b_{0}^{2}/R^{2}.
\label{923e12}
\end{equation}%
Thus Theorem \ref{p-harm} is established. $\Box $

\end{document}